\documentclass[12pt,reqno]{amsart}
\usepackage{amssymb,amsthm,url}
\usepackage[breaklinks=true]{hyperref}
\usepackage[utf8]{inputenc}
\usepackage[OT2,T1]{fontenc}
\swapnumbers
\let\<\langle
\let\>\rangle

\advance\hoffset by -0.5truein
\DeclareSymbolFont{cyrletters}{OT2}{wncyr}{m}{n}
\DeclareMathSymbol{\Sha}{\mathalpha}{cyrletters}{"58}

\theoremstyle{definition}
\newtheorem{definition}{Definition}[section]

\theoremstyle{plain}
\newtheorem{lemma}[definition]{Lemma}

\newtheorem{theorem}[definition]{Theorem}

\newenvironment{Proof}[1][Proof.]{\begin{trivlist}
\item[\hskip \labelsep {\bfseries #1}]}{\flushright
$\Box$\end{trivlist}}

\usepackage[all]{xy}
\usepackage{pdfsync}
\usepackage{ytableau}

\newcommand{\sudda}[1]{}

\begin{document}

\title{Simple $n$-Lie Poisson algebras}

\author{Farukh Mashurov}

\email{f.mashurov@gmail.com}
\address{SICM, Southern University of Science and Technology, Shenzhen, 518055, China}

\address{SDU University, Kaskelen, Kazakhstan}



\maketitle

\begin{abstract} 
Let $(A,\cdot,\omega)$ be a simple $n$-Lie Poisson algebra over a field of zero characteristic, $ 1 \in A.$ Then we prove that the $n$-Lie algebra $A^{[1]}/(A^{[1]}\cap Z)$ is simple, where $A^{[1]}$ denotes the derived $n$-Lie ideal and $Z$ is the center of $n$-Lie algebra $(A,\omega)$.
\end{abstract}

\section{\label{nn}\ Introduction}

In 1985, V.T. Filippov introduced $n$-Lie algebras as the generalization of Lie algebras \cite{Fil1985}. 
These algebras have attracted attention due to their relevance in the generalization of Hamiltonian mechanics proposed by Nambu \cite{Nambu}, later formalized mathematically by Takhtajan \cite{Takhtajan} as Nambu mechanics. Since then, $n$-Lie algebras have become an important area of interaction between mathematics and physics.

An \emph{$n$-Lie algebra} is a structure that generalizes the concept of a Lie algebra to an $n$-ary operation. Specifically, an $n$-Lie algebra over a field $F$ is a vector space $L$ equipped with an $n$-ary skew-symmetric bracket operation:
\[
\omega(\cdot, \cdot, \dots, \cdot) : L \times L \times \dots \times L \rightarrow L,
\]
such that for all $x_1, x_2, \dots, x_{n}, y_2, \dots, y_n \in L$, the following \emph{generalized Jacobi identity} holds:
$$\omega(\omega(x_1, x_2, \dots, x_{n}),y_2, \dots, y_{n})=\sum_{i=1}^n \omega(x_1, x_2, \dots, x_{i-1},\omega(x_i,y_2, \dots, y_{n}),x_{i+1},\ldots,x_n).$$

 Cohomological methods for $n$-Lie algebras were initiated by Takhtajan and later refined by Gautheron, Daletsky, and others \cite{DalTakh, Gautheron, Takhtajan, Takhtajan2}, who explored deformation theory and connections to Leibniz algebras. A rich theory of deformations and cohomology has since been developed, including formal 1-parameter deformations \cite{Makh16}, morphism cohomology \cite{AFM18}, and the structure induced by Lie algebras \cite{AKM11}.

Filippov’s foundational works  \cite{Fil1985, Fil1998} provided several key examples of $n$-Lie algebras, notably:

- the $(n+1)$-dimensional vector product algebra $A_{n+1}$;

- the Jacobian algebra defined on a commutative associative algebra $A$ with $n$ commuting derivations $D_1,\ldots,D_n$, and multiplication defined as follows:
$$Jac(u_1,u_2,\ldots,u_n)=\begin{array}{|ccc|}
  D_1 u_1   &  \cdots & D_1 u_n \\
    \cdots & \cdots & \cdots \\

      D_n u_1   &  \cdots & D_n u_n
\end{array},$$
for every $u_1,u_2,\ldots,u_n\in A.$

Another important infinite-dimensional example was given by Dzhumadil'daev \cite{DzhJac}, where only $n-1$ commuting derivations $D_1,\ldots, D_{n-1}$ are used, the $n$-ary bracket is determined by
$$W(u_1,u_2,\ldots,u_n)=\begin{array}{|ccc|}
    u_1   &  \cdots &  u_n \\
  D_1 u_1   &  \cdots & D_1 u_n \\
    \cdots & \cdots & \cdots \\
      D_{n-1} u_1   &  \cdots & D_{n-1} u_n
\end{array},$$ 
for every $u_1,u_2,\ldots,u_n\in A.$

The structural theory of finite-dimensional $n$-Lie algebras was developed in  \cite{Fil1998, Kas1987, Kas1991, Ling} and see references therein. In particular, Ling in \cite{Ling} proved that only one simple finite-dimensional $n$-Lie algebra exists over an algebraically closed field $F$ of characteristic $0.$ Finite and infinite-dimensional, irreducible, highest weight representations of the vector product algebra $A_{n+1}$ were studied in \cite{Balibanu, DzhRep}. 

Filippov and Pozhidaev introduced a class of monomial $n$-Lie algebras $A_G(f,t)$ \cite{Fil1998, Pozh98}, and a simplicity criterion was established in \cite{Pozh99}.   Cantarini and Kac \cite{CanKac2010, CanKac2016} classified all simple linearly compact $n$-Lie superalgebras and  (generalized) $n$-Lie-Poisson algebras over a field $F$ of characteristic 0. Their classification of simple $n$-Lie algebras includes four examples: one is the $n+1$-dimensional vector product $n$-Lie algebra, while the remaining three are infinite-dimensional $n$-Lie algebras, two of which were previously introduced by Filippov and Dzhumadil'daev. In the case of classification of simple generalized $n$-Lie-Poisson algebras, they proved that any simple linearly compact generalized $n$-Lie-Poisson algebra is equivalent either to the algebra introduced by Filippov or Dzhumadil'daev \cite{CanKac2016}. All irreducible continuous representations of these two infinite-dimensional simple linearly compact $n$-Lie algebras were classified in \cite{BM17} and \cite{BM19}.

Several constructions of $(n+1)$-Lie algebras from existing $n$-Lie structures have been developed \cite{BWLZ12, Omirov}.  An overview of various $n$-ary generalizations of Lie algebras, including $n$-Lie, $n$-Leibniz and generalized Poisson structures, as well as their applications in mathematical physics and cohomology properties, is presented in \cite{DeAI2010}. 
The simplicity of transposed Poisson $n$-Lie algebras was studied by the author in \cite{FM}.  
In particular, it was shown that a transposed Poisson $n$-Lie algebra $(A,\cdot,\omega)$ is simple if and only if $(A,\omega)$ is simple as an $n$-Lie algebra \cite{FM}.  

In this paper, we show that the analogue of the main theorem in \cite{AA19}  and \cite{H61} is also valid for $n$-Lie Poisson algebras.
\begin{theorem}\label{th1}
    Let $(A,\cdot,\omega)$ be a simple $n$-Lie Poisson algebra over a field of zero characteristic, $1\in A.$ Then the $n$-Lie algebra $A^{[1]}/(A^{[1]}\cap Z)$ is simple.
\end{theorem}

Throughout the paper, we assume that $(A,\cdot,\omega)$ is a simple $n$-Lie Poisson algebra over a field of characteristic~$0$ and that $1\in A$.

\section{Main results}

Let $A$ be an associative commutative algebra over a field $F$.  
An $n$-linear map
\[
\omega :\underbrace{A\times \cdots \times A}_n \longrightarrow A
\]
is called an \emph{$n$-ary Poisson bracket} if the pair $(A,\omega)$ forms an $n$-Lie algebra and the bracket satisfies the Leibniz rule
\begin{equation}\label{Leibniz rule}
    \omega(ab,u_2,\ldots,u_n)
    = a\,\omega(b,u_2,\ldots,u_n)
    + \omega(a,u_2,\ldots,u_n)\,b
\end{equation}
for all $a,b,u_2,\ldots,u_n\in A$.  
An algebra equipped with such a bracket is called an \emph{$n$-Lie Poisson algebra}.

In every $n$-Lie Poisson algebra the following identity holds:
\begin{equation}\label{Poisson ident}
    \omega(ab,c,u_3,\ldots,u_n)
    = \omega(a,bc,u_3,\ldots,u_n)
      + \omega(b,ac,u_3,\ldots,u_n).
\end{equation}

An ideal $I\subseteq A$ is called an \emph{$n$-ary Poisson ideal} if
\[
\omega(I,A,\ldots,A)\subseteq I.
\]
An $n$-Lie Poisson algebra that has no nontrivial $n$-ary Poisson ideals is called \emph{simple}.

We denote by
\[
Z=\{\,a\in A \mid \omega(a,A,\ldots,A)=0\,\}
\]
the center of the $n$-Lie algebra $(A,\omega)$.

For any $a_2,\ldots,a_n \in A$, consider the adjoint operator
\[
\operatorname{ad}(a_2,\ldots,a_n): A \longrightarrow A, \qquad
\operatorname{ad}(a_2,\ldots,a_n)(b)=\omega(b,a_2,\ldots,a_n).
\]
By the Leibniz rule \eqref{Leibniz rule}, each operator $ad(a_2,\ldots,a_n)$ is a derivation of the associative algebra $(A,\cdot)$.

Finally, define the derived series of ideals by
\[
A^{[0]} = A, \qquad
A^{[i+1]} = \omega(A^{[i]}, A^{[i]},\ldots,A^{[i]}), \quad i \ge 0.
\]

In what follows, we restate several lemmas established in \cite{AA19}. For clarity, we rewrite these results in the setting of $n$-ary multiplications.
 
\begin{lemma}\label{lem1}
Let $(A,\cdot,\omega)$ be a simple $n$-Lie Poisson algebra. Then the associative commutative algebra $A$ contains no nonzero nilpotent elements.
\end{lemma}
\begin{Proof}
Since the operator $ad$ is a derivation of the commutative associative algebra $A$.
Hence, the proof is the same as Lemma~1 in \cite{AA19}, and the claim follows.
\end{Proof}

\begin{lemma}\label{lem5}
Let $(A,\cdot, \omega)$ be a simple $n$-Lie Poisson algebra.  
Let $a_2,\ldots,a_n\in A$ and $m\ge 1$. If 
\[
ad(a_2,\ldots,a_n)^m = 0,
\]
then 
\[
ad(a_2,\ldots,a_n)=0.
\]
\end{lemma}
\begin{Proof}
  
The proof is identical to  Lemma~5 in \cite{AA19}. Let $b\in A$, and set
\[
b_1=\operatorname{ad}(a_2,\ldots,a_n)^{m-2}(b).
\] Indeed, applying the iterated derivation formula 
\[ ad(a_2,\ldots,a_n)^m(b_1^m)=
  m!\,\omega(b_1,a_2,\ldots,a_n)^m 
    + \]\[\sum_{\substack{k_1+\cdots+k_m=m \\ (k_1,\ldots,k_m)\ne(1,\ldots,1)}}
       {ad}(a_2,\ldots,a_n)^{k_1}(b_1)\cdots
       {ad}(a_2,\ldots,a_n)^{k_m}(b_1),
\]
and then using Lemma~\ref{lem1}, we obtain the desired result.

\end{Proof}

An ideal $U$ is called abelian if $\omega(U, U,\ldots,U)  =  (0).$

\begin{lemma}\label{lem6.0} Let $U$ be an abelian ideal of the algebra $A^{[1]}.$ Then $$\omega(A,  A^{[1]},\ldots,A^{[1]},U) = (0).$$ \end{lemma}

\begin{Proof} 
 We prove the statement by induction on the number of arguments $k$ from $A^{[1]}=\omega( A,\ldots,A )$. That is, for $0\leq k \leq n-2$, we show that
$$\omega(A,\underbrace{A^{[1]},\ldots,A^{[1]}}_{k},\underbrace{U,\ldots,U}_{n-k-1})=(0).$$

If $k=0,$ then for an arbitrary elements $a_2,\ldots,a_{n} \in U$ and  $a_1\in A$ we have
$$ad(a_2,\ldots,a_n)^3 (a_1) = \omega(\omega(\omega(a_1,a_2,\ldots,a_{n}),a_2,\ldots,a_{n}),a_2,\ldots,a_n) \in $$ 
$$\omega(\omega(\omega(A,\ldots A),U,U,\ldots,U),U,U,\ldots, U) \subseteq \omega(U, U, U\ldots,U)  =  (0).$$
By Lemma \ref{lem5} \begin{equation}\label{eq1}
ad(a_2,\ldots,a_n)(A) =\omega(A,a_2,\ldots,a_n)=(0).\end{equation}

Now, we show that $$\omega(A,A^{[1]},\underbrace{U,\ldots,U}_{n-2})=(0).$$
Again, for arbitrary elements $a_1\in A,$ $a_2\in A^{[1]},$  $a_3,\ldots,a_{n} \in U,$ and by \eqref{eq1}, we have
$$ad(a_2,a_3,\ldots,a_{n} )^3 (a_1) = \omega(\omega(\omega(a_1,a_2,\ldots,a_{n}),a_2,a_3\ldots,a_{n}),a_2,a_3,\ldots,a_n) \in $$ 
$$\omega(\omega(A^{[1]},A^{[1]},U,\ldots,U),A^{[1]},U,\ldots, U) \subseteq \omega(U, U\ldots,U,A)  =  (0).$$
By Lemma \ref{lem5} \begin{equation}\label{eq2} ad(a_2,a_3,\ldots,a_n)(A) =\omega(A,a_2,a_3,\ldots,a_n)=(0).\end{equation}

Suppose that  for all $k$ such that $1 \le k \le n-3$, we know
$$\omega(A,\underbrace{A^{[1]},\ldots,A^{[1]}}_{k},\underbrace{U,\ldots,U}_{n-k-1})=(0).$$
Then, we need to show that 
$$\omega(A,\underbrace{A^{[1]},\ldots,A^{[1]}}_{k+1},\underbrace{U,\ldots,U}_{n-k-2})=(0).$$
For $a_1\in A, a_2\ldots,a_{k+2} \in A^{[1]}$ and $a_{k+3},\ldots,a_{n}\in U$ we have
$$ad(a_2,\ldots,a_{n})^3 (a_1) = \omega(\omega(\omega(a_1,\ldots,a_{n}),a_2,\ldots,a_{n}),a_2,\ldots,a_{n}) \in$$ 
$$\omega(\omega(A^{[1]},\underbrace{A^{[1]},\ldots,A^{[1]}}_{k+1},\underbrace{U,\ldots,U}_{n-k-2}),\underbrace{A^{[1]},\ldots,A^{[1]}}_{k+1},\underbrace{U,\ldots,U}_{n-k-2})\subseteq$$
$$\omega(U,\underbrace{A^{[1]},\ldots,A^{[1]}}_{k+1},\underbrace{U,\ldots,U}_{n-k-2})\subseteq$$
$$\omega(A,\underbrace{A^{[1]},\ldots,A^{[1]}}_{k},\underbrace{U,\ldots,U}_{n-k-1})=(0).$$

By Lemma \ref{lem5} $$ad(a_2,\ldots,a_{n-1},u)(A) =\omega(A,A^{[1]},\ldots,A^{[1]},U)=(0),$$ which completes the proof of the lemma. 

\end{Proof}

\begin{lemma}\label{lem6} Let $U$ be an abelian ideal of the algebra $A^{[1]}.$ Then $$\omega( A,\ldots,A, U) = (0).$$ \end{lemma}

\begin{Proof}  Similarly to the previous lemma, we prove the statement by induction on the number of arguments $k$ from $A$ and using Lemma \ref{lem5} and Lemma \ref{lem6.0}. That is, for $1\leq k \leq n-1$, we show that
$$\omega(\underbrace{A,\ldots,A}_{k},\underbrace{A^{[1]},\ldots,A^{[1]}}_{n-k-1},U)=(0).$$
The case $k=1$ follows from Lemma \ref{lem6.0}.
Let us show that for $k=2.$

Let  $a_1,a_2\in A,a_3\ldots,a_{n-1}\in A^{[1]},$ $u \in U,$ and by Lemma \ref{lem6.0}, we have
$$ad(a_2,a_3,\ldots,a_{n-1},u)^2 (a_1) = \omega(\omega(a_1,a_2,a_3,\ldots,a_{n-1},u),a_2,a_3,\ldots,a_{n-1},u)\in $$ 
$$\omega(A,A^{[1]},A^{[1]},\ldots,A^{[1]},U)   =  (0).$$
By Lemma \ref{lem5} $$ad(a_2,a_3,\ldots,a_{n-1},u)(a_1) =\omega(a_1,a_2,a_3,\ldots,a_{n-1},u)=(0).$$ We have \begin{equation}\label{eq2} \omega(A,A,A^{[1]},\ldots,A^{[1]},U)=(0).\end{equation}

 Assume that  for all $k$ such that $1 \le k \le n-2$, we know that
$$\omega(\underbrace{A,\ldots,A}_{k},\underbrace{A^{[1]},\ldots,A^{[1]}}_{n-k-1},U)=0.$$
Then, we need to show that 
$$\omega(\underbrace{A,\ldots,A}_{k+1},\underbrace{A^{[1]},\ldots,A^{[1]}}_{n-k-2},U)=0.$$

For
$a_1,\ldots,a_{k+1}\in A,$ $ 
a_{k+2},\ldots,a_{n-1}\in A^{[1]}, u\in U$ and by Lemma \ref{lem6.0}, we have
$$ad(a_2,a_3,\ldots,a_{n-1},u)^2 (a_1) = \omega(\omega(a_1,a_2,a_3,\ldots,a_{n-1},u),a_2,a_3,\ldots,a_{n-1},u)\in $$ 
$$\omega(A^{[1]},\underbrace{A,\ldots,A}_{k},\underbrace{A^{[1]},\ldots,A^{[1]}}_{n-k-1},U)   =  (0).$$
By Lemma \ref{lem5} \begin{equation}\label{eq2} ad(a_2,a_3,\ldots,a_{n-1},u)(a_1) =\omega(a_1,a_2,a_3,\ldots,a_{n-1},u)=(0).\end{equation} Thus, $$\omega(\underbrace{A,\ldots,A}_{k+1},\underbrace{A^{[1]},\ldots,A^{[1]}}_{n-k-2},U)=0.$$

\end{Proof}

For a subspace $X \subset A$, we denote by
\[
{id}_{A}(X) = X + AX
\]
the ideal of $A$ generated by $X$.  
Let $U$ be an ideal of the Lie algebra $\omega(A,A,\ldots,A)$.  
We define the descending derived series of ideals by
\[
U^{[0]} = U, \qquad 
U^{[i+1]} = \omega(U^{[i]}, U^{[i]}, \ldots, U^{[i]}) \quad (i \ge 0).
\]

\begin{lemma}\label{lem7} 
    Let $U$ be an ideal of the algebra $A^{[1]}$ such that $U^{[3]}=(0).$ Then $$\omega(A,A,\ldots,A,U^{[1]})=(0).$$
\end{lemma}
\begin{Proof}
    The ideal $U^{[2]}$ of the $n$-Lie algebra $A^{[1]}$ is abelian. Hence by Lemma \ref{lem6} $$\omega(U^{[2]},A,\ldots,A)=(0).$$ 
    Let $u_2,\ldots,u_n\in U^{[1]},$

    $$ad(u_2,\ldots,u_n)^4(A)=\omega(\omega(\omega(\omega(A,u_2,\ldots,u_n),u_2,\ldots,u_n),u_2,\ldots,u_n),u_2,\ldots,u_n)$$
    $$\subseteq\omega(\omega(U^{[1]},\ldots,U^{[1]}),U^{[1]},\ldots,U^{[1]})\subseteq\omega(U^{[2]},U^{[1]},\cdots,U^{[1]})$$$$\subseteq\omega(U^{[2]},A,\cdots,A)=(0).$$

    By Lemma \ref{lem5} we have $$\omega(A,U^{[1]},\cdots,U^{[1]})=(0).$$

First, for $0\leq k \leq n-2$, we show that
$$\omega(A,\underbrace{A^{[1]},\ldots,A^{[1]}}_{k},\underbrace{U^{[1]},\ldots,U^{[1]}}_{n-k-1})=(0).$$

Let $k=1$ and  $a_1\in A,a_2\in A^{[1]}, a_3\ldots,a_{n}\in U^{[1]}$ and by Lemma \ref{lem6.0}, we have
$$ad(a_2,a_3,\ldots,a_{n})^3 (a_1) = $$
$$\omega(\omega(\omega(a_1,\ldots,a_n),a_2,\ldots,a_n),a_2,\ldots,a_n)\in $$ 
$$\omega(\omega(A^{[1]},A^{[1]},U^{[1]},\ldots,U^{[1]}),A^{[1]},U^{[1]},\ldots,U^{[1]})   \subseteq$$
$$\omega(U^{[1]},A^{[1]},U^{[1]},\ldots,U^{[1]})   =  (0).$$
By Lemma \ref{lem5} we have $$\omega(A,A^{[1]},U^{[1]},\ldots,U^{[1]})=(0).$$

Suppose that  for all $k$ such that $1 \le k \le n-3$, we know
$$\omega(A,\underbrace{A^{[1]},\ldots,A^{[1]}}_{k},\underbrace{U^{[1]},\ldots,U^{[1]}}_{n-k-1})=(0).$$
Then, we need to show that 
$$\omega(A,\underbrace{A^{[1]},\ldots,A^{[1]}}_{k+1},\underbrace{U^{[1]},\ldots,U^{[1]}}_{n-k-2})=(0).$$

For $a_1\in A, a_2\ldots,a_{k+2} \in A^{[1]}$ and $a_{k+3},\ldots,a_{n}\in U^{[1]}$ we have
$$ad(a_2,\ldots,a_{n})^3 (a_1) =$$$$ \omega(\omega(\omega(a_1,\ldots,a_{n}),a_2,\ldots,a_{n}),a_2,\ldots,a_{n}) \in$$ 
$$\omega(\omega(A^{[1]},\underbrace{A^{[1]},\ldots,A^{[1]}}_{k+1},\underbrace{U^{[1]},\ldots,U^{[1]}}_{n-k-2}),\underbrace{A^{[1]},\ldots,A^{[1]}}_{k+1},\underbrace{U^{[1]},\ldots,U^{[1]}}_{n-k-2})\subseteq$$
$$\omega(U^{[1]},\underbrace{A^{[1]},\ldots,A^{[1]}}_{k+1},\underbrace{U^{[1]},\ldots,U^{[1]}}_{n-k-2})\subseteq$$
$$\omega(A,\underbrace{A^{[1]},\ldots,A^{[1]}}_{k},\underbrace{U^{[1]},\ldots,U^{[1]}}_{n-k-1})=(0).$$

Thus, we proved that for $0\leq k \leq n-2$,
\begin{equation}\label{Aa1u1=0}   
\omega(A,\underbrace{A^{[1]},\ldots,A^{[1]}}_{k},\underbrace{U^{[1]},\ldots,U^{[1]}}_{n-k-1})=(0).
\end{equation}

Now, let  $a_1,a_2\in A,a_3\ldots,a_{n-1}\in A^{[1]},$ $a_{n} \in U^{[1]},$ and by \eqref{Aa1u1=0}, we have
$$ad(a_2,a_3,\ldots,a_{n-1},a_n)^2 (a_1) = \omega(\omega(a_1,a_2,a_3,\ldots,a_{n-1},a_n),a_2,a_3,\ldots,a_{n-1},u_n)\in $$ 
$$\omega(A^{[1]},A,A^{[1]},\ldots,A^{[1]},U^{[1]})   =  (0).$$
By Lemma \ref{lem5} \begin{equation}\label{eq2} ad(a_2,a_3,\ldots,a_{n-1},u_n)(a_1) =\omega(a_1,a_2,a_3,\ldots,a_{n-1},u_n)=(0).\end{equation} We have $$\omega(A,A,A^{[1]},\ldots,A^{[1]},U^{[1]})=(0).$$

 Assume that  for all $k$ such that $1 \le k \le n-3$, we know that
$$\omega(\underbrace{A,\ldots,A}_{k},\underbrace{A^{[1]},\ldots,A^{[1]}}_{n-k-1},U^{[1]})=0.$$
Then, we need to show that 
$$\omega(\underbrace{A,\ldots,A}_{k+1},\underbrace{A^{[1]},\ldots,A^{[1]}}_{n-k-2},U^{[1]})=0.$$

For
$a_1,\ldots,a_{k+1}\in A,$ $ 
a_{k+2},\ldots,a_{n-1}\in A^{[1]}, a_n\in U^{[1]},$  by hypothesis we have
$$ad(a_2,a_3,\ldots,a_{n-1},a_n)^2 (a_1) = \omega(\omega(a_1,a_2,a_3,\ldots,a_{n-1},u),a_2,a_3,\ldots,a_{n-1},u)\in $$ 
$$\omega(A^{[1]},\underbrace{A,\ldots,A}_{k},\underbrace{A^{[1]},\ldots,A^{[1]}}_{n-k-2},U^{[1]})   =  (0).$$
By Lemma \ref{lem5} $$ad(a_2,a_3,\ldots,a_{n-1},u_n)(a_1) =\omega(a_1,a_2,a_3,\ldots,a_{n-1},u_n)=(0).$$  Thus, $$\omega(\underbrace{A,\ldots,A}_{k+1},\underbrace{A^{[1]},\ldots,A^{[1]}}_{n-k-2},U^{[1]})=0.$$

\end{Proof}

\begin{lemma}\label{lem8} 
    Let $U$ be an ideal of the algebra $A^{[1]}$ such that $U^{[3]}=(0).$ Then $$\omega(A,\ldots,A,U)=(0).$$
\end{lemma}
\begin{Proof}

    Let $u_2,\ldots,u_n\in U,$ then by Lemma \ref{lem7}
    $$ad(u_2,\ldots,u_n)^4(A)=\omega(\omega(\omega(\omega(A,u_2,\ldots,u_n),u_2,\ldots,u_n),u_2,\ldots,u_n),u_2,\ldots,u_n)\subseteq$$    $$\omega(\omega(U,\ldots,U),U,\ldots,U)\subseteq\omega(U^{[1]},U,\cdots,U)\subseteq\omega(U^{[1]},A,\cdots,A)=(0).$$
    
     By Lemma \ref{lem5} we have $$\omega(A,U,\cdots,U)=(0).$$

Repeating the procedure in the proof of  Lemma~\ref{lem7}, we can obtain that 
$$\omega(\underbrace{A,\ldots,A}_{k},\underbrace{A^{[1]},\ldots,A^{[1]}}_{n-k-1},U)=0,$$
for every $0 \leq k \leq n-1.$

\end{Proof}

The following lemma is an analogue of Lemma 2 in \cite{AA19}.

\begin{lemma}\label{lem2}
    $\omega(id_A(U^{[3]}), A,\cdots,A) \subseteq U.$
\end{lemma}
\begin{Proof}
    For an arbitrary ideal $I$ of the $n$-Lie algebra $A^{[1]},$ by generalized Jacobi identity we have $$\omega(\omega(I,I,\cdots,I),A, \cdots,A)  \subseteq
\omega(\omega(I, A, \cdots,A),I,\cdots,I) \subseteq I.$$ 
 Hence,
$$\omega(U^{[i+1]}, A,\cdots,A) \subseteq U^{[i]}, \text{ for } i\geq 0.$$

By the identity \eqref{Poisson ident}
 \begin{equation}\label{ident4}
     \begin{array}{c}
     \omega(u\   \omega(a, v, u_3,\cdots,u_n),b,b_2,\cdots,b_n) = \\ 
    \omega(u,   \omega(a, v, u_3,\cdots,u_n)b,b_2,\cdots,b_n)+\omega(   \omega(a, v, u_3,\cdots,u_n),u b,b_2,\cdots,b_n)
 \end{array}
 \end{equation}

For $a,b,u, u_i,b_j\in A$ $(i=3,\ldots,n, j=2,3,\ldots,n),$ by \eqref{Leibniz rule} we have  $$
    \omega(\omega(au, v,u_3,\cdots,u_n) ,b,b_2,\cdots,b_n) = $$$$
      \omega(a\ \omega(u,v, u_3,\cdots,u_n),b,b_2,\cdots,b_n)  + \omega(u\   \omega(a, v, u_3,\cdots,u_n),b,b_2,\cdots,b_n).$$
     By applying identity \eqref{ident4} to the second term on the right-hand side, we obtain
$$ \omega(\omega(au, v,u_3,\cdots,u_n) ,b,b_2,\cdots,b_n) =  \omega(a\ \omega(u,v, u_3,\cdots,u_n),b,b_2,\cdots,b_n) + $$$$\omega(u,   \omega(a, v, u_3,\cdots,u_n)b,b_2,\cdots,b_n)+\omega(   \omega(a, v, u_3,\cdots,u_n),u b,b_2,\cdots,b_n).$$

Hence, $$
\begin{array}{c}
      \omega(a\ \omega(u,v, u_3,\cdots,u_n),b,b_2,\cdots,b_n) =\omega(\omega(au, v,u_3,\cdots,u_n) ,b,b_2,\cdots,b_n)-\\
     \omega(u,   \omega(a, v, u_3,\cdots,u_n)b,b_2,\cdots,b_n)-\omega(   \omega(a, v, u_3,\cdots,u_n),u b,b_2,\cdots,b_n).
\end{array}$$

The first and the third terms on the right-hand side lie in $\omega(\omega(v, A,\cdots, A), A,\cdots, A),$ the second summand lies in $\omega(u, A, A,\cdots, A)$. 

If $u, v, u_3,\ldots, u_n  \in U^{[2]}$ then by the remark above $$\omega(v, A,\cdots,A) + \omega(u,A, \cdots,A) \subseteq U^{[1]}$$ and $$\omega(\omega(v, A,\cdots,A),A,\cdots,A) \subseteq  \omega(U^{[1]}, A,\cdots,A)\subseteq U^{[0]} = U.$$ Therefore, $$ \omega(id_A(U^{[3]}), A,\cdots,A) \subseteq U,$$ which completes the proof of the lemma.
\end{Proof}

The following lemma is an $n$-ary analogue of Lemma 3 in \cite{AA19}, and its proof requires new techniques.

\begin{lemma}\label{lem3}
Let $I\subseteq A$ be a nonzero (associative commutative) ideal such that
\[
\omega(I,A^{[1]},\ldots,A^{[1]})\subseteq I.
\]
Then $I=A$.
\end{lemma}

\begin{Proof}  Let $\sqrt{I}$ be the set of all elements of the algebra $A$ that are nilpotent modulo $I.$ That is,
\[
\sqrt{I}
=
\{\,a\in A \mid a^{n}\in I \text{ for some } n\ge 1\,\}.
\]
This is an ideal of $A$, and $\sqrt{I}/I$ is the nilradical of the associative
commutative algebra $A/I$. 

Let $d$ be a derivation of $A$ such that $d(I)\subseteq I$.  
Then $d$ induces a derivation 
\[
\overline{d} : A/I \longrightarrow A/I.
\]
Since the nilradical of an algebra $A/I$ over a field of characteristic $0$ 
is invariant under all derivations  of $A/I$ \cite{Dix1997}, we have
\[
\overline{d}(\sqrt{I}/I)\subseteq \sqrt{I}/I,
\]
and therefore $d(\sqrt{I})\subseteq \sqrt{I}$.  
Taking $d=ad(x_1,\ldots,x_{n-1})$ for $x_1,\ldots,x_{n-1}\in\omega(A,A,\ldots,A)$, we obtain
\[
\omega(\sqrt{I},A^{[1]},\ldots,A^{[1]})\subseteq \sqrt{I}.
\]

If $\sqrt{I}=A$, then $1\in\sqrt{I}$, so $1$ is nilpotent modulo $I$.
Therefore $I=A$.
Thus, replacing $I$ by $\sqrt{I}$ if necessary, we may assume from now on that
\[
\sqrt{I}=I.
\]

For  $a\in A, a_3,\ldots,a_n \in A^{[1]},$ and $u \in I$ we have $$\omega(\omega(u, a^2,a_3,\ldots,a_n),u,a_3,\ldots,a_n) =$$
    $$-2\omega(u, a \  \omega(u, a,a_3,\ldots,a_n),a_3,\ldots,a_n)=$$
    $$-2(\omega(u,a,a_3,\ldots,a_n)^2+\omega(u, \omega(u, a,a_3,\ldots,a_n),a_3,\ldots,a_n)a).$$

    Hence $\omega(u,a,a_3,\ldots,a_n)^2 \in I,\omega(u,a,a_3,\ldots,a_n) \in \sqrt{I} = I.$ We have
$\omega(I , A,A^{[1]},\ldots,A^{[1]}) \subseteq I .$


  In the next step, for $a,a_3\in A, a_4,\ldots,a_n \in A^{[1]}, u \in I$ we have $$\omega(\omega(u, a^2,a_3,\ldots,a_n),u,a_3,\ldots,a_n) =$$ $$-2\omega(u, a \  \omega(u, a,a_3,\ldots,a_n),a_3,\ldots,a_n)=$$ $$-2(\omega(u,a,a_3,\ldots,a_n)^2+\omega(u, \omega(u, a,a_3,\ldots,a_n),a_3,\ldots,a_n)a).$$


  Hence $\omega(u,a,a_3,\ldots,a_n)^2 \in I,\omega(u,a,a_3,\ldots,a_n) \in \sqrt{I} = I,$ and $$\omega(A,A,\underbrace{A^{[1]},\ldots,A^{[1]}}_{n-3},I)\subseteq I.$$

  Proceeding inductively, one finally obtains
$\omega(I , A,\ldots,A) \subseteq I ,$ i.e $I$ is a $n$-ary Poisson ideal of $A.$ Since the $n$-Lie Poisson algebra $(A, \cdot, \omega )$ is simple, we conclude that $I = A$. This completes the proof of the lemma.

\end{Proof}

\begin{Proof}[The proof of Theorem \ref{th1}] Assume $U$ is a proper ideal of the $n$-Lie algebra $A^{[1]}.$
Let $I = id_A(U^{[3]}).$ It is clear that $$\omega(I,A^{[1]},A^{[1]},\ldots,A^{[1]}) \subseteq I.$$ If $U^{[3]}\neq 0$ then by Lemma \ref{lem3} we have $I = A.$ Hence by Lemma \ref{lem2}  $$\omega(A, A,\ldots, A)\subseteq U.$$ Therefore, we have that $U^{[3]}=0.$

Finally, by Lemma~\ref{lem8}, every ideal $U$ of $A^{[1]}$ such that $U^{[3]}=0$ is contained in the center $Z$. 
\end{Proof}

\section*{Acknowledgements}

The author sincerely thanks Professor Efim Zelmanov for posing the problem that motivated this work. This research was carried out during the author’s postdoctoral stay at the Southern University of Science and Technology (SUSTech). The author is also grateful to the referee for valuable comments and suggestions. 

\end{document}